\documentclass[sn-mathphys]{sn-jnl}% Math and Physical Sciences Reference Style
\usepackage[utf8]{inputenc}
\usepackage[english]{babel}
\usepackage[T1]{fontenc}
\usepackage{amsmath}
\usepackage{bm}
\usepackage{siunitx}
\usepackage{graphicx}

\usepackage{setspace}

\usepackage{cancel}
\usepackage{balance}

\usepackage{listings}

\usepackage{placeins}
\usepackage{float}
\usepackage{subfig}
\usepackage{natbib}
\usepackage{epstopdf}
\graphicspath{{graph/}{data/}}
\definecolor{color0}{HTML}{E6001A}
\definecolor{color1}{HTML}{0083CC}
\definecolor{color2}{HTML}{EE7A34}
\definecolor{color3}{HTML}{50B695}

\jyear{2022}%

\theoremstyle{thmstyleone}%
\theoremstyle{thmstyletwo}%

\theoremstyle{thmstylethree}%

\begin{document}
\title[Transient Forward Harmonic Adjoint Sensitivity Analysis]{Transient Forward Harmonic Adjoint Sensitivity Analysis}
	
%	%%=============================================================%%
%	%% Prefix	-> \pfx{Dr}
%	%% GivenName	-> \fnm{Joergen W.}
%	%% Particle	-> \spfx{van der} -> surname prefix
%	%% FamilyName	-> \sur{Ploeg}
%	%% Suffix	-> \sfx{IV}
%	%% NatureName	-> \tanm{Poet Laureate} -> Title after name
%	%% Degrees	-> \dgr{MSc, PhD}
%	%% \author*[1,2]{\pfx{Dr} \fnm{Joergen W.} \spfx{van der} \sur{Ploeg} \sfx{IV} \tanm{Poet Laureate} 
%	%%                 \dgr{MSc, PhD}}\email{iauthor@gmail.com}
%	%%=============================================================%%
%	
	\author*[1]{\fnm{Julian} \sur{Sarpe}}\email{julian\_johannes.buschbaum@tu-darmstadt.de}
	\author[2]{\fnm{Andreas} \sur{Klaedtke}}\email{andreas.klaedtke@de.bosch.com}
%	\equalcont{These authors contributed equally to this work.}
%	
	\author[1]{\fnm{Herbert} \sur{De Gersem}}\email{degersem@temf.tu-darmstadt.de}
%	\equalcont{These authors contributed equally to this work.}
%	
	\affil*[1]{\orgdiv{Institute for Accelerator Science and Electromagnetic Fields (TEMF)}, \orgname{Technische Universit\"at Darmstadt}, \orgaddress{\street{Schlo{\ss}gartenstra{\ss}e~8}, \city{64289 Darmstadt}, \country{Germany}}}
	\affil[2]{\orgdiv{Corporate Sector Research and Advance Engineering}, \orgname{Robert Bosch GmbH}, \orgaddress{\street{Robert-Bosch-Campus 1}, \city{71272 Renningen}, \country{Germany}}}

	\abstract{%
		This paper presents a transient forward harmonic adjoint sensitivity analysis (TFHA), which is a combination of a transient forward circuit analysis with a harmonic balance based adjoint sensitivity analysis. 
		TFHA provides sensitivities of quantities of interest from time-periodic problems w.r.t. many design parameters, as used in the design process of power-electronics devices.
		The TFHA shows advantages in applications where the harmonic balance based adjoint sensitivity analysis or finite difference approaches for sensitivity analysis perform poorly. 
		In contrast to existing methods, the TFHA can be used in combination with arbitrary forward solvers, i.e. general transient solvers. 
	}
	
	\maketitle
	
	\section{Introduction}
		The development cost of devices containing electronic circuits is directly linked to the number of prototyping cycles.
		A means of reducing the number of prototype cycles is the use of electric circuit analysis tools, such as Spice, that help predict the behavior of the circuits without the need to create any hardware. 
		For time-harmonic systems, the analysis can be performed either in time or frequency domain. 
		Transient circuit analysis treats nonlinear devices by linearizing their behavior at every time step. 
		In frequency domain, nonlinear devices cause a mutual dependence of different spectral components.
		Thus, frequency domain analyses are not generally possible~\cite{nakhla1976piecewise, kundert1986nonlin}.  
		Nonetheless, a nonlinear circuit problem can be solved in frequency domain using the harmonic balance (HB) method~\cite{nakhla1976piecewise}. 
		The HB method approximates the solution of the steady state for a finite number of harmonics with Newton iteration~\cite{maas2003nonlinear}. 
		The HB method avoids the issue of long transients in time domain simulations.
		However one also has to consider the drawbacks of HB.
		The more harmonics that are necessary to approximate the nonlinear device behavior, the more degrees of freedom (DoF) the equation system has. 
		Moreover, the convergence of the Newton procedure is suffering in strongly nonlinear systems~\cite{rizzoli1983general, maas2003nonlinear}. 
		\par
		A systematic approach for circuit analysis is the employment of
		sensitivity analysis, which aims to reduce optimization cycles in the early stages of development by systematically optimizing certain design parameters~\cite{nikolova2004adjoint}. 
		A common global sensitivity analysis is performed by computing Sobol indices from a Polynomial Chaos Expansion (PCE) surrogate model.~\cite{SOBOL2001271}. 
		However, PCE and similar methods are limited by the number of design parameters due to the curse of dimensionality~\cite{SUDRET2008964} that makes this type of analysis unfeasible for large parameter spaces.
		Another category of sensitivity analysis techniques involves computing the gradient of the Quantity of Interest (QoI) w.r.t. one or more design parameters. 
		Gradient based methods are also referred to as local sensitivity analysis~\cite{SUDRET2008964}. 
		\par
		Adjoint sensitivity analysis is the method of choice for a setting with many design parameters~\cite{cao2003Adjoint, nikolova2004adjoint}. 
		Transient adjoint sensitivity analysis~\cite{cao2003Adjoint, ilievski2007adjoint} as well as sensitivity analysis based on HB solvers~\cite{bandler1988unified} have been published and used before. 
		Transient adjoint sensitivity analysis is an efficient approach for the analysis of strongly nonlinear systems.
		However, when dealing with time-dependent sensitivities, the inherent part of solving multiple adjoint problems creates a computational bottleneck~\cite{nikolova2004adjoint}. 
		This issue does not exist for HB based sensitivity analysis. 
		Combining the advantages of transient analysis and HB based adjoint sensitivity analysis, it is possible to utilize an HB based adjoint sensitivity analysis while at the same time avoiding convergence problems and many necessary Newton iterations for strongly nonlinear problems. 
		\par
		This paper proposes the transient forward adjoint harmonic sensitivity analysis (TFHA) as a combination of a nonlinear transient forward analysis with a harmonic adjoint sensitivity analysis. 
		The TFHA is an efficient and robust method for sensitivity analysis in power-electronics applications as demonstrated in practical examples in section~\ref{sec:results}. 
		\par
		This paper is structured as follows: Section~\ref{sec:circuit} gives a brief introduction to the used circuit analysis methods and necessary considerations for the analysis of nonlinear circuits.
		Section~\ref{sec:sens} gives an overview over the mentioned sensitivity analysis methods. Section~\ref{sec:TFHA} introduces the TFHA as a novel sensitivity analysis method.
		After that, TFHA is applied to several model problems in section~\ref{sec:results}. 
		Conclusions on the applicability of the TFHA are drawn in section~\ref{sec:conclusions}.
	
%	\clearpage
	\section{Nonlinear Circuit Analysis}
		\label{sec:circuit}
		To perform the sensitivity analysis, a simulation framework for the considered nonlinear circuits must be introduced first. 
		All presented methods are based on the modified nodal analysis (MNA), which will be introduced first.
		\subsection{Modified Nodal Analysis}
			MNA has been the method of choice for most circuit problems since its introduction in 1975~\cite{ruehli1975mna}.
			MNA extends nodal analysis to accommodate impedance devices such as inductances or voltage sources in addition to the admittance devices.
			In the static case, the MNA problem reads:
			\begin{equation}
				\bm{F}(\bm{x}) = \bm{A} \bm{x} - \bm{i}_\mathrm{s} = 0.
				\label{eq:statMNA}
			\end{equation}
			The solution vector $\bm{x}$ contains the nodal voltages $\bm{u}_\mathrm{nodes}$ for all nodes and the edge currents $\bm{i}_\mathrm{edge}$ for the impedance devices.
			MNA employs an extended system matrix $\bm{A}$, containing the nodal admittance matrix $\bm{Y}$, the edge impedance matrix $\bm{Z}$ of the impedance devices and two incidence matrices $\bm{B}$ and $\bm{B}^\mathsf{T}$ coupling both submatrices:	
			\begin{equation}
				 \bm{A} = \left(
				 \begin{array}{c|c}
				 \begin{matrix}
				 & &  & & \\
				 & & \bm{Y} &  \\
				 & & & & 
				 \end{matrix}
				 & 
				 \begin{matrix}
				 \\
				 \bm{B}^\mathsf{T}\\
				 \\
				 \end{matrix}
				 \\\hline
				 \begin{matrix}
				 & & \bm{B} & &
				 \end{matrix}	
				 &
				 \begin{matrix}
				 \bm{Z}
				 \end{matrix}
				 \end{array}
				 \right).
				\label{eq:mnamat}
			\end{equation}
			The right hand side $\bm{i}_\mathrm{s}$ of the nodal analysis is modified accordingly such that it contains all input nodal currents $\bm{i}_\mathrm{nodes}$ and edge voltages $\bm{u}_\mathrm{edge}$ for the respective impedance devices:
			\begin{equation}
				\bm{i}_\mathrm{s} = \left(
				\begin{array}{c|c}
				\begin{matrix} & & \bm{i}_\mathrm{nodes} & &  \end{matrix} & \bm{u}_\mathrm{edge}
				\end{array}
				\right)^\mathsf{T}.
				\label{eq:mnarhs}
			\end{equation}

		\subsection{Nonlinear Transient Circuit Analysis}
			MNA can be utilized for a transient circuit analysis.
			The transient MNA is given as a differential algebraic equation (DAE) system 
			\begin{equation}
				\bm{F}(\bm{x},\dot{\bm{x}},t) =
		 		\bm{A}_\mathrm{C} \dot{\bm{x}}(t) + \bm{A}_\mathrm{G} \bm{x}(t) - \bm{i}_\mathrm{s}(t) = 0, \mathrm{~~~~~} \bm{x}(t=0) = 0,
				\label{eq:TDMNA}
			\end{equation}
			where $\bm{A}_\mathrm{C}$ contains capacitor and inductor contributions, $\bm{A}_\mathrm{G}$ contains conductance contributions and $\bm{i}_\mathrm{s}$ contains the independent sources.
			If the circuit contains nonlinear devices, a linearization of the system is required. 
			This is done by using the Newton method~\cite{vlach1983computer} and the introduction of the Jacobian matrix $\bm{J}$
			\begin{equation}
				\bm{J}_\mathrm{G}(\bm{x}, t) = \bm{A}_\mathrm{G} - \frac{\partial \bm{i}_\mathrm{nl}(\bm{x}(t), t)}{\partial \bm{x}(t)},
				\label{eq:jacobian}
			\end{equation} 
			where $\bm{i}_\mathrm{nl}$ denotes the vector of voltage dependent currents for nonlinear admittance devices.
			For impedance and energy storage devices, the linearization is performed analogously. 
			For transient analyses, the linearization is performed for each time step $t$.
	
		\subsection{Harmonic Balance Method}
			A linear system can directly be transferred to frequency domain~\cite{nikolova2004adjoint}: 
			\begin{equation}
				\underline{\bm{F}}(\underline{\bm{x}},\omega) = \mathrm{j}\omega \bm{A}_\mathrm{C} \underline{\bm{x}}(\omega) + \bm{A}_\mathrm{G} \underline{\bm{x}}(\omega) - \underline{\bm{i}}_\mathrm{s}(\omega) = 0,
				\label{eq:ACMNAraw}
			\end{equation}
			where the underlining indicates the phasors of the corresponding quantities and $\omega$ is the circular frequency. 
			The resulting system matrices are combined into a single matrix $\underline{\bm{A}} = \mathrm{j}\omega\bm{A}_\mathrm{C} + \bm{A}_\mathrm{G} $ which gives the consitutive equation for the default MNA in frequency domain:
			\begin{equation}
				\underline{\bm{F}}(\underline{\bm{x}}) = \underline{\bm{A}} \underline{\bm{x}} - \underline{\bm{i}}_\mathrm{s} = 0.
				\label{eq:ACMNA}
			\end{equation}
			The HB method is a formalism to approximate the solution of a nonlinear system iteratively with a Newton iteration in frequency domain~\cite{nakhla1976piecewise}.
			The Jacobian~\eqref{eq:jacobianAC} is defined at multiple frequencies in order to represent the nonlinear behavior
			\begin{equation}
				\underline{\bm{J}}(\bm{x}, \bm{\omega}) = \underline{\bm{A}}(\bm{\omega}) - \frac{\partial \underline{\bm{i}}_\mathrm{nl}(\underline{\bm{x}}, \bm{\omega})}{\partial \underline{\bm{x}}},
				\label{eq:jacobianAC}
			\end{equation}
			where $\bm{\omega}$ is the vector of all considered frequencies in the system.
			Approximation with multiple harmonic frequencies increases the number of degrees of freedom (DoF) the stronger the nonlinearities~\cite{maas2003nonlinear}.
			Resultingly, the HB method is very performant in weakly nonlinear systems that exhibit long transient times.
			This performance advantage is diminished for strongly nonlinear systems due to the exploding number of DoFs.
			\par
			The holomorphy of the function vector~\eqref{eq:ACMNA} needs to be ensured in the computational domain to define the Jacobian matrix. 
			The Paley-Wiener theorem~\cite{hille1933annals} states that the Fourier transform of a bounded function $F\in L^2$ is holomorphic if the function variable is in the upper half plane, i.e. when t is restricted to $\Re_+$,
			which can be assumed in real world applications without the loss of generality.
			
	\section{Sensitivity Analysis}
		This section summarizes known methods to obtain the sensitivity for a QoI $U$ w.r.t. one or multiple design parameters $p$.
		The QoI can be any quantity within the circuit such as an edge voltage or current. 
		Design parameters in our application are components within the circuit, such as inductances, capacitances or resistances.
		The presented methods are limited to gradient based sensitivities.
		This limitation is introduced since the applications exhibit very large parameter spaces which makes global sensitivity analysis not well suited~\cite{SUDRET2008964}. 
		\label{sec:sens}
		\subsection{Transient Direct Sensitivity Analysis}
			$U(\bm{x})$ depends on the circuit solution, and hence, indirectly on time.
			Therefore, the calculation of the sensitivity~\cite{cao2003Adjoint}:
			\begin{equation}
				\frac{\mathrm{d} U}{\mathrm{d} p} = \frac{\partial U}{\partial \bm{x}} \frac{\mathrm{d} \bm{x}}{\mathrm{d} p}
			\end{equation}
			involves the calculation of the derivative $\mathrm{d} \bm{x}/\mathrm{d} p$ of the circuit solution $\bm{x}$ w.r.t. the parameter $p$. 
			$\partial U/\partial \bm{x}$ is a mapping operator that obtains the sensitivity for the QoI from the derivative $\frac{\mathrm{d} \bm{x}}{\mathrm{d} p}$.
			$\frac{\mathrm{d} \bm{x}}{\mathrm{d} p}$ is obtained by taking the derivative of~\eqref{eq:TDMNA} w.r.t. $p$. 
			\begin{equation}
				\frac{\mathrm{d} \bm{F}}{\mathrm{d} p} = \bm{A}_\mathrm{C} \frac{\mathrm{d} \dot{\bm{x}}}{\mathrm{d} p} + \bm{J}_\mathrm{G} \frac{\mathrm{d} \bm{x}}{\mathrm{d} p} + \left(\frac{\mathrm{d} \bm{A}_\mathrm{C}}{\mathrm{d} p}\dot{\bm{x}} + \frac{\mathrm{d} \bm{A}_\mathrm{G}}{\mathrm{d} p}\bm{x}\right) = 0.
				\label{eq:theprob}
			\end{equation}
			This approach is referred to as direct sensitivity analysis (DSA)~\cite{nikolova2004adjoint}. 
			Eq.\eqref{eq:theprob} needs to be solved individually for every design parameter. 
			As a result, DSA is expensive for the sensitivity analysis w.r.t. many design parameters. 

			\subsection{Transient Adjoint Sensitivity Analysis}
\begin{figure}[b]
				\centering
				\includegraphics[width=\linewidth]{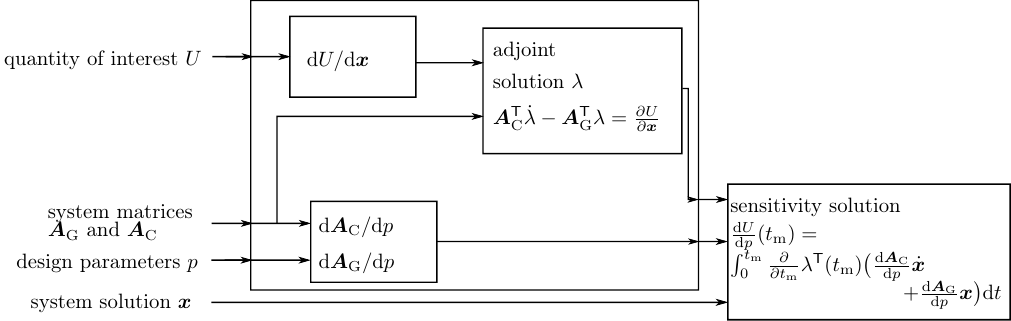}
				\caption{Workflow of the transient adjoint sensitivity method.}
				\label{fig:flowadjoinharmTD}
			\end{figure}

			Transient adjoint sensitivity analysis determines sensitivities based on the system solution $\bm{x}$ and a test function $\bm{\lambda}$. 
			The workflow for transient adjoint sensitivity analysis is outlined in Fig.~\ref{fig:flowadjoinharmTD}. 
%			\par
			To obtain the adjoint sensitivity analysis, differential equation~\eqref{eq:theprob} is multiplied with $\lambda$ as a test function and integrated over time:
			\begin{equation}
				\int_0^{t_\mathrm{m}} \lambda^\mathsf{T} (t_\mathrm{m}) \left(\bm{A}_\mathrm{C} \frac{\mathrm{d} \dot{\bm{x}}}{\mathrm{d} p} + \bm{J}_\mathrm{G} \frac{\mathrm{d} \bm{x}}{\mathrm{d} p}\right)~\mathrm{d} t = - \int_0^{t_\mathrm{m}}  \lambda^\mathsf{T}(t_\mathrm{m}) \left(\frac{\mathrm{d} \bm{A}_\mathrm{C}}{\mathrm{d} p}\dot{\bm{x}} + \frac{\mathrm{d} \bm{A}_\mathrm{G}}{\mathrm{d} p}\bm{x}\right) ~\mathrm{d} t .
			\end{equation}
			Integration by parts eliminates the time derivative of $\mathrm{d} \bm{x}/\mathrm{d} p$~\cite{nikolova2004adjoint}: 
			\begin{multline}
				\int_0^{t_\mathrm{m}} \left( - \dot{\lambda}^\mathsf{T} \bm{A}_\mathrm{C}  + \lambda^\mathsf{T}(t_\mathrm{m}) \bm{J}_\mathrm{G}  \right)\frac{\mathrm{d} \bm{x}}{\mathrm{d} p} ~\mathrm{d} t
				\\= - \int_0^{t_\mathrm{m}} \lambda^\mathsf{T}(t_\mathrm{m}) \left(\frac{\mathrm{d} \bm{A}_\mathrm{C}}{\mathrm{d} p}\dot{\bm{x}} + \frac{\mathrm{d} \bm{A}_\mathrm{G}}{\mathrm{d} p}\bm{x}\right) ~\mathrm{d} t - \left[ \lambda^\mathsf{T}(t_\mathrm{m}) \bm{A}_\mathrm{C} \frac{\mathrm{d} \bm{x}}{\mathrm{d} p}  \right]_0^{t_\mathrm{m}}.
				\label{eq:between}
			\end{multline}
	 		Without loss of generality $\lambda$ is chosen as the solution of the adjoint system, to eliminate the boundary terms and accommodate the left hand side term:
			\begin{equation}
				\bm{A}_\mathrm{C}^\mathsf{T} \dot{\lambda} - \bm{J}_\mathrm{G}^\mathsf{T} \lambda = \frac{\partial U}{\partial \bm{x}}, \mathrm{~~~~~} \lambda (t=t_\mathrm{m}) = 0.
				\label{eq:adjoint}
			\end{equation}
		 	The adjoint solution is calculated backwards in time to ensure the reverse initial condition $\lambda (t=t_\mathrm{m}) = 0$ for arbitrary $t_\mathrm{m}$~\cite{cao2003Adjoint}.
			The term at $t=0$ is eliminated due to the initial condition $\bm{x}(t=0) = 0$ in Eq.\eqref{eq:TDMNA}.
			Substituting Eq.\eqref{eq:adjoint} into Eq.\eqref{eq:between} and multiplication by -1 yields the integration for the sensitivity of the QoI $U$.
			\begin{equation}
				\int_0^{t_\mathrm{m}} \frac{\mathrm{d} U}{\mathrm{d} p} ~\mathrm{d} t = \int_0^{t_\mathrm{m}} \lambda^\mathsf{T}(t_\mathrm{m}) \left(\frac{\mathrm{d} \bm{A}_\mathrm{C}}{\mathrm{d} p}\dot{\bm{x}} + \frac{\mathrm{d} \bm{A}_\mathrm{G}}{\mathrm{d} p}\bm{x}\right)~\mathrm{d} t
				\label{eq:sensintegral}
			\end{equation}
			The generalized version of Leibniz integral rule is applied to the right hand side term~\cite{Flanders1973Differentiation}, to obtain the sensitivity for the QoI at a specific time instant $t_\mathrm{m}$:
			\begin{equation}
				\frac{\mathrm{d} U}{\mathrm{d} p} (t_\mathrm{m}) = \int_0^{t_\mathrm{m}} \frac{\partial}{\partial T} \lambda^\mathsf{T}(t_\mathrm{m}) \left(\frac{\mathrm{d} \bm{A}_\mathrm{C}}{\mathrm{d} p}\dot{\bm{x}} + \frac{\mathrm{d} \bm{A}_\mathrm{G}}{\mathrm{d} p}\bm{x}\right)~\mathrm{d}t.
			\end{equation}
			There is an adjoint solution for each time instant $t_\mathrm{m}$, as computed by solving equation~\eqref{eq:adjoint}.
			Thus, the task of computing the solution to Eq.\eqref{eq:adjoint} has to be redone for each instant in time where the sensitivity is to be considered.
			
	\subsection{Harmonic Balance based Direct Sensitivity Analysis}
		\label{sec:sensharm}
		Analogously to the transient case, DSA is derived by symbolic differentiation of the reverse initial condition from Eq.\eqref{eq:ACMNA} for the HB method. 
		\begin{equation}
			\frac{\mathrm{d}\underline{\bm{F}}(\bm{\varphi})}{\mathrm{d} p} = \frac{\mathrm{d}\underline{\bm{A}}}{\mathrm{d} p} \underline{\bm{x}} + \left(\underline{\bm{A}}  - \frac{\partial\underline{\bm{i}}_\mathrm{nl}}{\partial\underline{\bm{x}}}\right) \frac{\mathrm{d}\underline{\bm{x}}}{\mathrm{d} p} = \frac{\mathrm{d}\underline{\bm{A}}}{\mathrm{d} p} \underline{\bm{x}} + \underline{\bm{J}} \frac{\mathrm{d}\underline{\bm{x}}}{\mathrm{d} p} = 0
			\label{eq:senswithjac}
		\end{equation}
		From Eq.\eqref{eq:senswithjac} we can solve the system for the sensitivity $\mathrm{d}\underline{\bm{x}}/\mathrm{d}p$:
		\begin{equation}
			\frac{\mathrm{d}\underline{\bm{A}}}{\mathrm{d} p} \underline{\bm{x}} + \underline{\bm{J}} \frac{\mathrm{d}\underline{\bm{x}}}{\mathrm{d} p} = 0 ~~~~~\Rightarrow~~~~~  \frac{\mathrm{d}\underline{\bm{x}}}{\mathrm{d} p} = - \underline{\bm{J}}^{-1} \frac{\mathrm{d}\underline{\bm{A}}}{\mathrm{d} p} \underline{\bm{x}}
			\label{eq;directMatMult}
		\end{equation}
		The Jacobian is determined in each Newton iteration when used for the harmonic balance solver. 
		If the goal is to perform sensitivity analysis based on an existing circuit solution, the Jacobian is given as the converged solution from the nonlinear problem~\cite{bandler1988unified}. 
		\par
		If many design parameters $p$ are analyzed, the numerical efficiency of HB based direct sensitivity decreases, because Eq.\eqref{eq;directMatMult} contains a matrix multiplication with the dense matrix $\underline{\bm{J}}^{-1}$.
		This dense matrix multiplication can be avoided by introducing an HB based adjoint sensitivity analysis.
		
		\subsection{Harmonic Balance based Adjoint Sensitivity Analysis}
		\label{sec:HBadj}
		The adjoint system for the harmonic balance case~\eqref{eq:senswithjac} is given as
		\begin{equation}
			\underline{\bm{J}}^\mathsf{H} \bm{\underline{\lambda}} = \frac{\partial \underline{U}}{\partial\underline{\bm{x}}} ~~~~~\Rightarrow~~~~~  \bm{\underline{\lambda}} = \frac{\partial \underline{U}}{\partial\underline{\bm{x}}} \underline{\bm{J}}^{-\mathsf{H}}.
			\label{eq:harmadj}
		\end{equation}
		As the Jacobian does not depend on the adjoint solution $\bm{\lambda}$, the adjoint problem is linear.
		The sensitivity of the QoI $\underline{U}(p)$ w.r.t. the design parameter $p$ follows as:
		\begin{equation}
			\frac{\mathrm{d} \underline{U}(p)}{\mathrm{d} p} = -\left(\frac{\partial \underline{U}}{\partial\underline{\bm{x}}}\right)^\mathsf{H} \underline{\bm{J}}^{-1} \frac{\mathrm{d}\underline{\bm{A}}}{\mathrm{d} p} \underline{\bm{x}} = - \bm{\underline{\lambda}}^\mathsf{H} \frac{\mathrm{d}\underline{\bm{A}}}{\mathrm{d} p} \underline{\bm{x}}.
			\label{eq:adjointSensHarm}
		\end{equation}
		In frequency domain, the adjoint solution can be reused for all considered frequencies, which eliminates the issues of the transient analysis but increases the number of DoFs for large frequency spaces.
		Compared to the HB based direct sensitivity analysis, HB based adjoint sensitivity analysis does not contain any dense matrix multiplications, improving performance.
		The procedure is outlined as the green dashed box in Fig.~\ref{fig:flowadjoinharm}.
		\par
		The issue that remains is the bad performance for strongly nonlinear systems.
		This is particularly problematic when the Jacobian has to be approximated along many convergence steps of the HB iteration (Eq.\eqref{eq:ACMNA}).
		A hybrid time-frequency domain method, that eliminates the issue of many Newton iteration steps, is proposed in the following section.

	\section{Transient Forward Harmonic Adjoint Sensitivity Analysis}
	\label{sec:TFHA}
		The HB method is well suited for simulating steady-state operation.
		If strong nonlinearities occur, the harmonics become coupled across a wide range of frequencies.
		Especially in the context of power electronics, short rise and fall times are present, which requires us to consider a large number of harmonics~\cite{maas2003nonlinear}. 
		\begin{figure}[b]
			\centering
			\includegraphics[width=\linewidth]{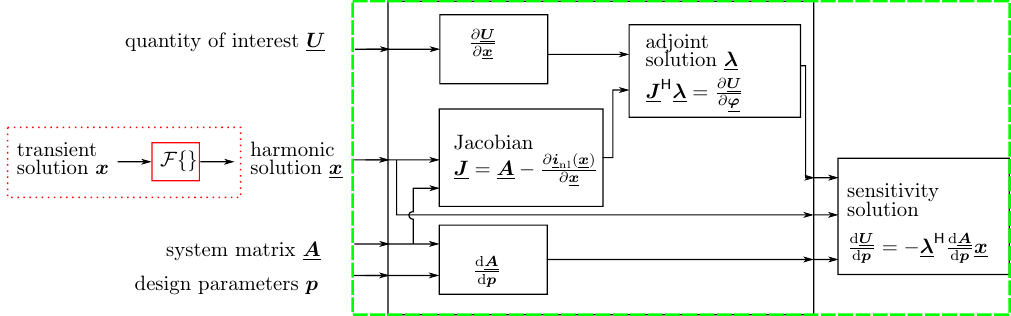}
			\caption{Workflow for the harmonic balance based adjoint sensitivity analysis and the TFHA.}
			\label{fig:flowadjoinharm}
		\end{figure}
		A procedure consisting of a forward HB solve combined with a harmonic adjoint method for calculating  the sensitivities was proposed in~\cite{bandler1988unified} and has been recapitulated in section~\ref{sec:sensharm}.
		Obtaining the forward solution for $\bm{x}$ with the harmonic balance method is not compulsory for the sensitivity analysis. 
		This motivates the following procedure. 
		The circuit solution is calculated by a transient circuit solver. 
		The adjoint solution is obtained by a harmonic analysis. 
		This procedure extends the HB based adjoint sensitivity analysis (in green) as depicted by the red dotted box in Fig.~\ref{fig:flowadjoinharm}.
		\par 
		The TFHA obtains the circuit solution $\bm{x}$ through the utilization of an efficient transient solver of choice. 
		Once an approximation to the steady state is found, the Fourier transform $\bm{\underline{x}}$ of one period of $\bm{x}$ is used to calculate the Jacobian $\bm{\underline{J}}$. 
		The approximated Jacobian matrix is then used to obtain the adjoint solution in frequency domain, analogously to the way presented in section~\ref{sec:HBadj}. 
		Subsequently, the adjoint solution is combined with the symbolic derivative of system matrix $\mathrm{d}\underline{\bm{A}}/\mathrm{d} p$ to calculate the sensitivities.
		Depending on the application, the sensitivity can be illustrated as a frequency spectrum or be transformed back to time.
		\par
		The TFHA is exposed to information loss  which results in a residual error, due to a finite number of harmonics. 
		The Euclidean distance of the solution with fewer harmonics against the solution with more harmonics quantifies the residual.
		This is based on the idea of the Zienkiewicz-Zhu error estimator~\cite{zienkiewicz1992superconvergent}.
		The relative error is approximated by the quotient of the absolute error divided by the norm of the fine solution:
		\begin{equation}
			\mathcal{E}_\mathrm{rel} \approx \frac{\left\| \left(\frac{\mathrm{d} \underline{U}}{\mathrm{d} p}\right)_\mathrm{fine} - \left(\frac{\mathrm{d} \underline{U}}{\mathrm{d} p}\right)_\mathrm{coarse} \right\|_{2}}{\left\| \left(\frac{\mathrm{d} \underline{U}}{\mathrm{d} p}\right)_\mathrm{fine}\right\|_{2}}.
			\label{eq:err_rel}
		\end{equation}
		Without prior considerations of the spectral circuit behavior, the sensitivity is determined by iteratively increasing the number of harmonics.
		The quantification of the error (Eq.\eqref{eq:err_rel}) is finally used as a termination condition to assess the converged solution.
			
	\section{Results}
		\label{sec:results}
		The TFHA is illustrated for three different examples:
		A half-wave rectifier as a nonlinear, academic example and
		two different power electronics circuits serving as industrial examples.
		\subsection{Half-Wave Rectifier}
\begin{figure}[b]
				\centering
				\includegraphics[width=.3\textwidth]{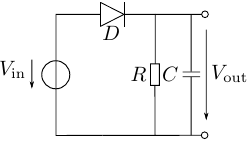}
				\caption{Schematic of a functional half-wave rectifying circuit.}
				\label{fig:HWrectcirc}
			\end{figure}

\begin{figure}[b]
				\centering
				\subfloat[Time domain plot over two periods]{\includegraphics[width=0.45\textwidth]{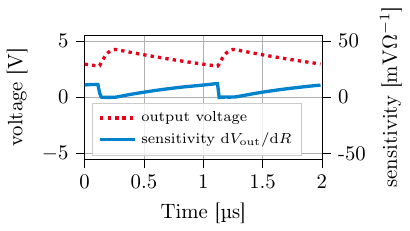}}
				\hspace{.04\textwidth}
				\subfloat[Frequency domain plot of sensitivity]{\includegraphics[width=0.45\textwidth]{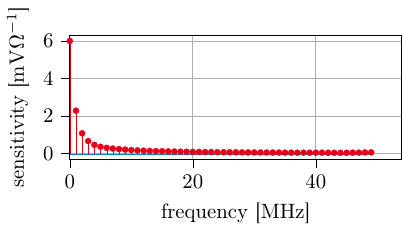}}
				\caption{Sensitivity of the output voltage  w.r.t. the resistance $R$ for the half-wave rectifying circuit in time and frequency domain.}
				\label{fig:TDhwsens}
			\end{figure}

			A half-wave rectifier rectifier (Fig.~\ref{fig:HWrectcirc}) consists of only one nonlinear and two linear components.
			Fig.~\ref{fig:TDhwsens} displays both the time domain and the frequency domain sensitivity of the output voltage w.r.t. the resistance $R$.
			Due to the weak nonlinearities within the circuit, the sensitivity solution converges fastly even for a small number of harmonics.
			
		\subsection{Boost Converter}
\begin{figure}[t]
				\centering
				\includegraphics[width=.65\linewidth]{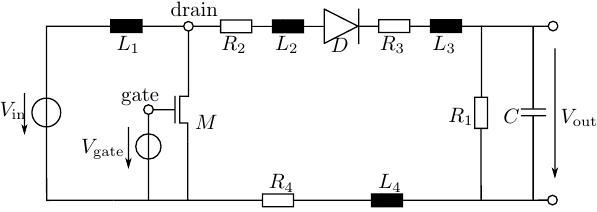}
				\caption{Boost converter circuit model with parasitics. The elements $L_{1}$, $R_{1}$, $C$ and $M$ model the functional behavior of the boost converter, the elements $R_{2}$, $L_{2}$, $R_{3}$, $L_{3}$, $R_{4}$ and $L_{4}$ are parasitic elements which model the EMC effects on the circuit board.}
				\label{fig:boostConvPar}
			\end{figure}

\begin{figure}[t]
				\centering
				\includegraphics[width=0.8\textwidth]{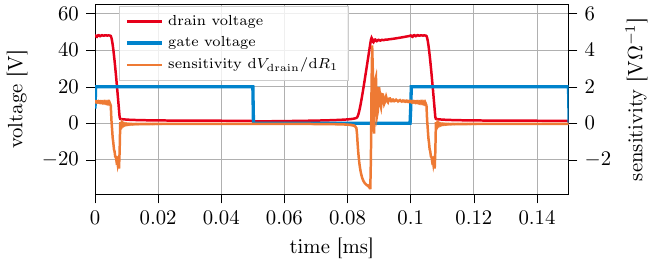}
				\caption{Time resolved sensitivity of the gate voltage $V_\mathrm{drain}$ w.r.t. the resistance $R_1$ related to the drain and gate voltages of the boost converter MOSFET.}
				\label{fig:boostParSens}
			\end{figure}

\begin{figure}[b]
				\centering
				\includegraphics[width=0.8\textwidth]{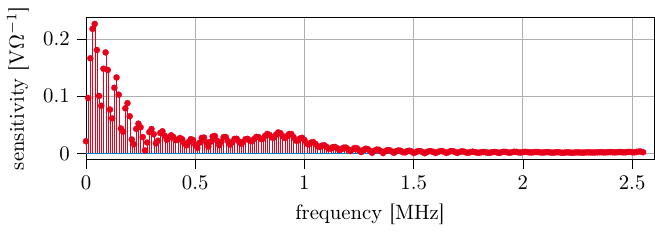}
				\caption{Sensitivity spectrum $\mathrm{d} V_\mathrm{drain} / \mathrm{d} R_1$ of the boost converter circuit showing a significant spectral content at around $\SI{0.8}{\mega\hertz}$ that must not be neglected.%
				 }
				\label{fig:boostR1spec}
			\end{figure}

			A boost converter is a type of DC-DC converter circuit which is used in many power electronic applications.
			In a circuit model, parasitic effects can be modeled as lumped circuit elements such as resistances, inductances and capacitances.
			The circuit with parasitic lumped elements is shown in Fig.~\ref{fig:boostConvPar}.
			The sensitivity of the drain voltage $V_\mathrm{drain}$ w.r.t. the output resistance $R_1$ is shown in Fig.~\ref{fig:boostParSens}.
			A significant overshoot is observed for the sensitivity and the drain voltage throughout the switching process.
			This is attributed to the broadband nature of the system.
			The overshoot can heavily influence the sensitivities especially close to the switching of the transistor. 
			Hence, this is a relevant test for accuracy and performance of the algorithm. The results in Fig.~\ref{fig:boostParSens} show that the algorithm correctly approximates the overshoot while limiting the Gibbs ripples when a frequency domain solution with an insufficient number of harmonics is calculated. 
%			which can distort the solution.
			The spectral TFHA solution gives a good indication of this behavior (Fig.~\ref{fig:boostR1spec}).
			The spectrum does not decrease steadily for higher harmonics, but rather shows a side lobe, with higher peaks at around $\SI{0.8}{\mega\hertz}$.
			 
		\subsection{Active Filter Circuit}
\begin{figure}[b]
				\centering
				\includegraphics[width=.75\linewidth]{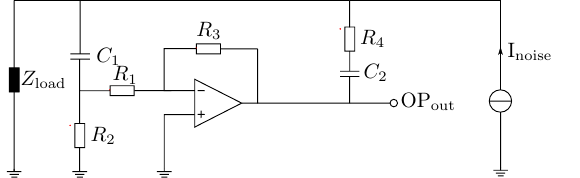}
				\caption{Functional schematic for the nonlinear active filter circuit.}
				\label{fig:filterNetwork}
			\end{figure}

\begin{figure}[b]
				\centering
				\includegraphics[width=\linewidth]{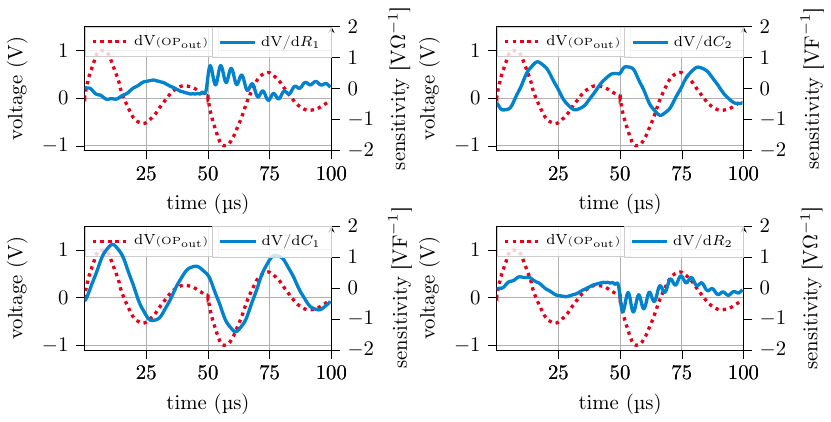}
				\caption{Sensitivity results for V(OP$_\mathrm{out}$) of the nonlinear active filter circuit w.r.t. different design parameters.}
				\label{fig:plotsFilter}
			\end{figure}

			An active filter circuit that injects a reversed polarity noise signal serves as an example with a larger number of elements.
			This circuit is nonlinear.
			Fig.~\ref{fig:filterNetwork} shows a reduced schematic of the entire active filter setup as used in the analysis.
			The complete circuit consists of 240 linear and nonlinear elements and has 152 degrees of freedom.
			The circuit contains several specific transistor and diode models. 
			The Xyce~\cite{Keiter2022Xyce} circuit simulator is used to assemble the Jacobian in frequency domain.
			The sensitivity of several QoIs is analyzed w.r.t. several design parameters of the circuit. The first QoI is the output potential of the OPAMP V(OP$_\mathrm{out})$.
			Due to the filter properties of the circuit, the quantity OP$_\mathrm{out}$ is smoother than other quantities in the circuit (Fig.~\ref{fig:filterNetwork}).
			The sensitivity results for the QoI V(OP$_\mathrm{out}$) w.r.t.the capacitances $C_1$ and $C_2$ as well as the resistances $R_1$ and $R_2$ are shown in Fig.~\ref{fig:plotsFilter}.
			The capacitances show a large influence on the amplitude of V($\mathrm{OP}_\mathrm{out}$). The resistances strongly influence the pre-filter properties on the input of the OPAMP. As a result, variations of the resistances lead to sensitivities which vary quickly in time. It is conclusive that the resistances have the most influence on the flank steepness of the OPAMP output.
			\par
			The sensitivity analysis is executed for all design parameters in the circuit to benchmark the performance of the TFHA. 
			The adjoint solution is calculated in 20 seconds using the implementation of Eq.\eqref{eq:harmadj}. 
			The evaluation of the sensitivity with Eq.\eqref{eq:adjointSensHarm} takes less than 2 seconds for each design parameter. 
			This is a significant speedup compared to existing methods for the given example. Direct sensitivity analysis with Xyce takes 12 seconds for each design parameter.
			As a result, a speedup of  $\SI{85}{\%}$ can be observed for the considered number of design parameters. 
			The sensitivities for the ten most influential design parameters are listed in table~\ref{tab:sens}. Note here that the design parameters not visible in Fig.~\ref{fig:filterNetwork} are part of the load impedance $Z_\mathrm{load}$.
			\begin{table}[t]
				\centering
				\caption{List of the ten most influential design parameters (ten largest normalized sensitivities).}
				\label{tab:sens}
				\begin{tabular}{| l | l |}
					\hline
					sensitivity & value (normalized to maximum sensitivity) \\ \hline \hline
					$\mathrm{d}V(OP_\mathrm{out}) / \mathrm{d} L_4$ & 1 \\ \hline
					$\mathrm{d}V(OP_\mathrm{out}) / \mathrm{d} C_1$ & 0.9733 \\ \hline
					$\mathrm{d}V(OP_\mathrm{out}) / \mathrm{d} R_3$ & 0.9590 \\ \hline
					$\mathrm{d}V(OP_\mathrm{out}) / \mathrm{d} C_2$ & 0.4331 \\ \hline
					$\mathrm{d}V(OP_\mathrm{out}) / \mathrm{d} C_6$ & 0.3670 \\ \hline
					$\mathrm{d}V(OP_\mathrm{out}) / \mathrm{d} C_3$ & 0.3664 \\ \hline
					$\mathrm{d}V(OP_\mathrm{out}) / \mathrm{d} R_{11}$ & 0.2990 \\ \hline
					$\mathrm{d}V(OP_\mathrm{out}) / \mathrm{d} C_5$ & 0.1264 \\ \hline
					$\mathrm{d}V(OP_\mathrm{out}) / \mathrm{d} R_2$ & 0.1056 \\ \hline
					$\mathrm{d}V(OP_\mathrm{out}) / \mathrm{d} R_1$ & 0.0982 \\ \hline
				\end{tabular}
			\end{table}
%			\clearpage
			\begin{figure}[t]
				\centering
				\includegraphics[width=.9\linewidth]{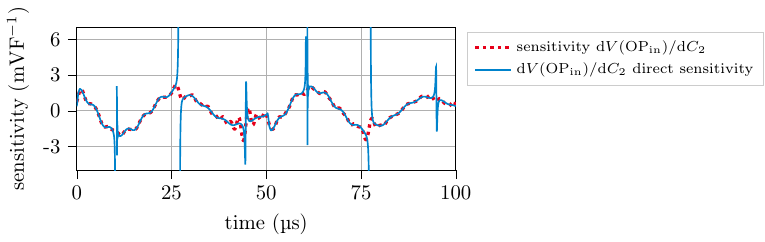}
				\caption{Sensitivity results for V($\mathrm{OP}_\mathrm{in}$) w.r.t. capacitance $C_2$ of the filter circuit computed by the TFHA compared to the transient direct Xyce~\cite{osti_1562422} sensitivity solution.}
				\label{fig:plotsFilterInput}
			\end{figure}
			The input node of the active filter has more oscillations at high frequencies. 
			As Fig.~\ref{fig:plotsFilterInput} shows, the V($\mathrm{OP}_\mathrm{in}$) sensitivity w.r.t. $C_2$ requires more spectral components for a high frequency approximation.
			To illustrate this, the sensitivity of V($\mathrm{OP}_\mathrm{in}$) w.r.t. element $C_2$ is calculated with the TFHA and compared to the sensitivity solution of the direct sensitivity solver of Xyce~\cite{osti_1562422}. 
			The comparison is shown in Fig.~\ref{fig:plotsFilterInput}. 
			As expected, TFHA can accurately determine the sensitivity in smooth intervals but deviates for sharp peaks.

	\FloatBarrier
	\section{Conclusions}
	\label{sec:conclusions}
		This paper recapitulated transient adjoint sensitivity analysis methods as well as HB based adjoint sensitivity analysis for nonlinear circuits. 
		Transient forward harmonic adjoint sensitivity analysis (TFHA) which combines transient solvers with harmonic adjoint sensitivity analysis was introduced as a novel method that improves the efficiency of sensitivity computation with regards to a large number of parameters. 
		Results calculated by the TFHA have been presented and compared to results calculated by transient or harmonic adjoint sensitivity analysis. 
		\par
		TFHA is advantageous in applications where purely transient sensitivity analysis is costly, but the HB method itself might need many Newton iterations to converge. 
		This can be the case in fastly oscillating or strongly nonlinear systems~\cite{maas2003nonlinear}, such as power electronic applications. 
		The TFHA requires one steady state time domain solution obtained by a solver of choice and only one solution for the adjoint system and therefore one inversion of the Jacobian matrix. 
		The TFHA performs well in weakly nonlinear systems where many design parameters are considered. 
		Additionally, TFHA is particularly performant for sensitivity analysis on time-dependent QoIs, since it requires only one adjoint solution in such scenarios.

	\section*{Declarations}
	\paragraph{Ethical Approval} Not applicable.
	\paragraph{Competing interests} The authors declare no conflict of interest.
	\paragraph{Authors' contributions} J. Sarpe wrote the main manuscript text and prepared the figures. A. Klaedtke provided the models for the "Boost Converter" and the "Active Filter Circuit". All authors wrote the introduction and the conclusions. All authors reviewed the manuscript.
	\paragraph{Funding} The authors declare no funding.
	\paragraph{Availability of data and materials} Not applicable.
%	
%	%************* Verzeichnisse ************
	\bibliography{refs}% common bib file

%% 	\clearpage
%	\printbibliography %Literatur
%	%\listoffigures %Abbildungen
%	%****************************************
%
%	%****************************************
%	% Anhang
%	%****************************************	
%% 	\appendix
%% 	\addcontentsline{toc}{chapter}{Anhang}

\end{document}